\documentclass[10pt]{amsart}

\usepackage{amsmath}
\usepackage{amssymb}
\usepackage{epsfig}
\usepackage{psfrag}

\newcommand\bbR{\mathbb{R}}
\newcommand\bbC{\mathbb{C}}

\newcommand{\beq}{\begin{equation}}
\newcommand{\eeq}{\end{equation}}
\newcommand{\beqn}{\begin{eqnarray}}
\newcommand{\eeqn}{\end{eqnarray}}
\newcommand{\tn}{\textnormal}

\newcommand{\p}{\psi}
\newcommand{\eps}{\varepsilon}

\newcommand{\Om}{\Omega}

\newcommand\nh{{\mathcal L(H^2)}}
\newcommand\nl{{\mathcal L(L^2)}}

\newtheorem{lemma}{Lemma}

\newtheorem{theorem}[lemma]{Theorem}

\newtheorem{rem}{Remark}
\newtheorem{algo}{Algorithm}

\title{Analysis of the Toolkit method for the time-dependant Schr\"odinger equation}

\begin{document}

\maketitle

\author{Lucie Baudouin,Julien Salomon,Gabriel Turinici}

\begin{abstract}
The goal of this paper is to provide an analysis of the ``toolkit'' method used in the
numerical approximation of the time-dependent Schr\"odinger equation. The ``toolkit'' method
is based on precomputation of elementary propagators and was seen to be very efficient in the optimal
control framework. Our analysis shows that this
method provides better results than the second order Strang
operator splitting. In addition, we present two improvements of the method in
the limit of low and large intensity control fields.    
\end{abstract}




\thispagestyle{plain}
\section{Introduction}

The control of the evolution of molecular systems at the quantum level
has been a long standing goal ever since the beginning of the laser technology.
After an initial slowed down of the investigations in this area due to unsuccessful experiments,
the realization
that the problem can be recast and attacked with the tools of
(optimal) control theory~\cite{hbref110} greatly contributed to the
first positive experimental
results~\cite{hbref1,hbref3,hbref10,hbref9,hbref4}. 
Ever since, the desire to understand theoretically how the laser acts to
control the molecule lead the investigators to resort to numerical simulations which
require repeated resolution of the Time Dependent Schr\"odinger Equation of the type~\eqref{schrod_original}; additional
motivation comes from related contexts (online
identification algorithms, learning algorithms, quantum computing~\cite{hbref90}, etc.).

The numerical method used to solve the time dependent Schr\"odinger equation must
provide accurate results without prohibitive computational cost. The conservation of the $L^2$ norm of the wave function $\p(x,t)$ 
is also generally required for stability and as a mean of qualitative validation
of the numerical solution.

In this context, the second order Strang operator splitting is often
considered~\cite{descombes,handbook9,strangsplitting}. However, this method suffers from two drawbacks. First, the numerical error is proportional to the norm of the control which
implies poor accuracy when dealing with large laser fields $\eps(t)$ and
make necessary the use of small time steps.
Secondly, it requires at each time step three matrix
products. This difficulty is enhanced in some particular settings e.g., in optimal control,
where the matrices involved in the control term must be assembled online.

Recently introduced, the ``toolkit method''~\cite{yip-mazziotti-rabitz-03,yip-mazziotti-rabitz-03bis} solves this last problem by
precomputing a set of elementary matrices, used in the numerical
resolution. Each matrix is associated to (one or several) field values and enables to solve the evolution over one time 
step.
This algorithm has been used in various frameworks and
shows excellent results. It has also been coupled successfully with optimal control and
identification issues~\cite{BST}. The dependence on the $L^\infty$-norm of the control, which is
a restriction of the Strang method, is also improved by the ``toolkit method'' as it will be shown in our analysis.

The goal of the paper is to provide a (first) numerical analysis of the ``toolkit method''. 
Our mathematical tools are related to that in~\cite{descombes} (but for a different setting; see also \cite{lubich,sanzserna}
for connected results); the treatment here is different 
because of the quantization appearing in the values of the control $\eps(t)$ 
which impacts both the mathematical analysis and
the numerical efficiency of the method.
The analysis enables us to propose two possible improvements. 

The paper is organized as follows: 
after having introduced the model and some notations in
Section~\ref{sec:modnot}, the toolkit method is presented and analyzed in 
Section~\ref{sec:tk}. An improvement of this method in the limit of small
control fields is introduced in Section~\ref{sec:itk}. A second
improvement, in the limit of large control fields is given in
Section~\ref{sec:itk2}. Finally, Section~\ref{sec:nr} gathers some numerical
results.

\section{Model and notations}\label{sec:modnot}

In this section, we present the Schr\"odinger Equation that will be
considered in the paper and some useful notations.

We consider the time dependent Schr\"odinger equation (TDSE)
\begin{eqnarray}\label{schrod_original}
\left\lbrace
\begin{array}{ll}
i\partial_t \p(x,t) = (H_0(x) - \mu(x)\eps(t)) \p(x,t),&\bbR^3\times(0,T)\\
\p(x,0)=\p_0(x),&\bbR^3.
\end{array} \right.
\end{eqnarray}
This equation governs the evolution of a quantum system, described by
its wave function $\p$, that interacts with a laser pulse of amplitude $\eps$, the control variable. The
factor $\mu$ is the dipole moment operator of the system. The Hamiltonian of the system is $H_0=-\Delta_x + V$ where $\Delta_x$ is
the Laplacian operator over the space variables and $V=V(x)$ the electrostatic potential in which
the system evolves. We refer to \cite{rabitz3} for more details about
models involved in quantum control.
Note that to obtain Eq. \eqref{schrod_original}, one has considered
the laser effect as a perturbative term, so that
the control term
$\eps (t)\mu(x)$ is obtained
through a first order approximation with respect to $\eps (t)$. While
often considered, this
approximation fails at describing some models
involving non linear laser-dipole interaction, see
e.g. \cite{claudedionjcp}. Consequently, the norm
of the field cannot be always considered as a small parameter, and
numerical solvers have to tolerate large controls, as the one described here after. \\

\indent Throughout this paper, $T>0$ is  the time of
control of a quantum system. The space $L^{p}(0,T;X)$, with $p\in 
[1,+\infty)$ denotes the usual Lebesgue space taking its values in a Banach
space $X$. The notation $W^{1,1}(0,T)$ corresponds to the space of time dependent functions belonging to $L^1(0,T;\bbR)$ such that their first time derivative also belongs to $L^1(0,T)$.
We denote by $L^2$ the space $L^2(\bbR^3,\bbC)$ and by $W^{2,\infty}$
and $H^2$ the Sobolev spaces $W^{2,\infty}(\bbR^3,\bbR)$ and
$H^2(\bbR^3,\bbC)$. The space $\mathcal L(H^2)$ is the space of linear
functionals on $H^2$. One can refer to \cite{Brezis} (or the introduction of \cite{caze}) for more details about the definitions of these functional spaces.\\

Finally, in order to introduce some numerical solver of \eqref{schrod_original}, let us
consider  an integer $N$ and $\Delta t>0$ such that
$N\Delta t=T$. We introduce the time discretization $(t_j)_{0\leq j\leq N}$ of
$[0,T]$ with $t_j=j\Delta t$ and we also denote by $t_{j+\frac 12}$ the intermediate time $\frac{t_j + t_{j+1}}2 = (j+ \frac 12)\Delta t$. \\

Let us first recall some basic results of existence and regularity of the solution of the TDSE.
These are corollaries of a general result on time dependent Hamiltonians (see \cite{R-S}, p.~285,
Theorem~X.70). \\
   
\begin{lemma}\label{schrod0}
Let $\mu \in \mathcal L(H^2)$, $V\in W^{2,\infty}$, $\eps\in L^{2}(0,T)$ and $\psi_0\in H^2$.
The Schr\"odinger equation
\beq\label{S}
\left\{
\begin{array}{ll}
i\partial_t \p (t) =\left(H_0-\mu\eps(t)\right)\p (t)&\bbR^3\times(0,T)\\
\p(0)=\p_0,&\bbR^3,
\end{array}
\right.
\eeq
has a unique solution $\p\in L^\infty(0,T;H^2)\cap W^{1,\infty}(0,T;L^2)$ such that 
\beq\nonumber
\|\p (t)\|_{L^\infty(0,T;H^2)} +  \|\partial_t \p
(t)\|_{L^\infty(0,T;L^2)} \leq C\|\mu\|_\nh\|\eps\|_{W^{1,1}(0,T)}
\|\psi_0\|_{H^2}.\eeq 
Moreover, for all $t\in [0,T],\quad  \|\p (t)\|_{L^2}=\|\psi_0\|_{L^2}$.\\
\end{lemma}

It is also well known (see \cite{caze} for instance) that for any
$T>0$ and $\phi_0 \in H^2$, if we have $\eps(t) = \bar\eps \in \bbR $,
independent of time $t$, the Schr\"odinger equation 
$$
\left\lbrace
\begin{array}{ll}
i\partial_t \phi (t) =\left(H_0-\mu\bar\eps\right)\phi (t),&\bbR^3\times(0,T)\\
\phi(0)=\phi_0,&\bbR^3
\end{array} \right.
$$
has a unique solution $\phi(t)=S(t)\phi_0$ such that $\phi\in
C([0,T];H^2)\cap C^1([0,T];L^2)$, where $(S(t))_{t\in \mathbb R}$
denotes the one-parameter semi-group generated by the operator
$H_0-\mu\bar\eps$. Moreover, for all $t\in [0,T]$, $S(t) \in \mathcal
L (H^2)$ and  we have  
\begin{equation}
\label{SSS}
\begin{array}{ll}
S(t)\phi_0 \in C(0,T;H^2),& \forall \phi_0 \in H^2;\\
\|S(t)\|_{ \mathcal L (H^2)} \leq 1+CT \leq K,& \forall t\in[0,T]; K=K(\|\mu\|_\nh,\eps_{\tn{max}});\\
\|S(t)\phi_0\|_{L^2} = \|\phi_0\|_{L^2},& \forall \phi_0 \in L^2, \forall t\in\bbR;\\
S(0) = \tn{Id}, \\
S(t+s) = S(t)S(s),& \forall s,t\in \bbR.\\
\end{array}
\end{equation}
Therefore, the solution of Eq.~\eqref{S} is
obtained equivalently as a solution to the integral equation 
$$ \p(t)=S(t) \p_0 + i \int_0^t S(t-s)(\eps(s)-\bar\eps)\mu \p(s)\,ds.$$

\section{The toolkit method}\label{sec:tk}

We now present the toolkit method and describe the corresponding error analysis. 

\subsection{Algorithm}

In this method, we assume that the control field $\eps$ satisfies the following
hypothesis:

\begin{equation}\tag{$\mathcal H$}
\label{H}
\forall t \in [0,T],\ \eps(t)\in [\eps_{\min},\eps_{\max}].
\end{equation}
The values of the control field are discretized 
according to: 
\beq\label{tkdisc}
\bar{\eps}_\ell=\eps_{\min}+\ell\Delta \eps,\ \ell=0...m,
\eeq 
with $m=\frac{\eps_{\max}-\eps_{\min}}{\Delta \eps}$.
Here, the values $\bar{\eps}_\ell$ have been here uniformly chosen in the
interval $[\eps_{\min},\eps_{\max}]$. If some properties of the field are
known, e.g. its mean value or its variance, some improvement of the
method can be obtained by optimizing the distribution 
of the values $\bar{\eps}_\ell$. More generally, this topic enters the
field of scalar quantization, that will not be considered
in this paper. We refer to \cite{sayood} and the references
therein for a review of standard methods in this domain.
In order to solve numerically equation~\eqref{S},
the toolkit method proceeds as follows.\\
\newpage
\begin{algo} (toolkit method)\label{tk}
\begin{enumerate}
\item Preprocessing. Precompute the ``toolkit'', i.e. the set of propagators:
$$S_{\ell}(\Delta t) \textnormal{ for } \ell = 0,\cdots, m,$$
where $(S_{\ell}(t))_{t\in \mathbb R}$ denotes the one-parameter
semi-group generated by the operator $H_0-\mu\bar{\eps}_\ell$, the
sequence $(\eps_\ell)_{\ell=0,\cdots, m,}$ being defined by
\eqref{tkdisc}. 
\item Given a control field $\eps\in L^2$ satisfying \ref{H} and
  $\psi^{K}_0=\psi_0$, the sequence $(\psi^{K}_{j})_{j=0,...,N}$ that
  approximates $(\psi(t_j))_{j=0,...,N}$, is obtained recursively by
  iterating the following loop:
\begin{enumerate} 
\item \label{step:mp}Find:
$$\ell_j = {\rm
 argmin}_{\ell=1,\cdots,m}\{|\eps(t_{j+\frac 12})-\bar{\eps}_\ell|\},$$
\item Set $\psi^{K}_{j+1}=S_{\ell_j}(\Delta t) \psi^{K}_j$.\\
\end{enumerate} 
\end{enumerate}
\end{algo}

In this ``toolkit'' approximation, we consider that the changes in the
Hamiltonian $H(t):=H_0-\mu\eps(t)$ can be neglected over a time step $\Delta
t$. In this way, if $\Delta \eps=0$  (infinite toolkit), and for
a relevant time discretization, the simulation corresponding to
piecewise constant control fields is exact. Such a property does not hold with 
methods that approximate the exponential, e.g. the second order Strang
operator splitting. Indeed, these approaches introduce an algebraic error, due to the 
non-commutation of the operators $H_0$ and $\mu$ that is consequently
proportional to $\|\eps\|_{L^\infty(0,T)}$.

\begin{rem}
In the original form of the toolkit method \tn{\cite{kurti-manby-ren-artamonov-ho-rabitz-05,
 yip-mazziotti-rabitz-03, yip-mazziotti-rabitz-03bis}}, the mid-point
choice proposed in Step~\ref{step:mp} of Algorithm~\ref{tk} is not considered. Yet, the
introduction of this strategy enables us to improve the order of the
method (see the analysis hereafter). 
\end{rem}

\subsection{Analysis of the method}

Let us now present an error analysis of the toolkit
method. More precisely, this section aims at proving the following result:\\
 
\begin{theorem}\label{anatk}
Let $\eps\in W^{2,\infty}(0,T)$ and $\psi$ the corresponding solution of
\eqref{S}. Let $\psi^{K}$ be the approximation of $\psi$ obtained
with Algorithm \ref{tk}. Given $\Delta t>0$ and $\Delta \eps>0$,
there exists $\lambda_1>0$, $\lambda_2>0$, that do not depend on $\|\eps\|_{L^\infty(0,T)}$ such that: 
\beq\label{est1}
\| \psi(T)-\psi^{K}(T)\|_{L^2}\leq \lambda_1 \Delta \eps  + \lambda_2\Delta t^2.
\eeq
Moreover, there exists $\nu_1>0$, $\nu_2>0$ depending on $\|\eps\|_{W^{1,1}(0,T)}$ such that: 
\beq\label{est2}
\| \psi(T)-\psi^{K}(T)\|_{H^2}\leq \nu_1 \Delta \eps  + \nu_2\Delta t^2.\\
\eeq
\end{theorem}
\begin{rem}
This result shows that the toolkit method enables to work with large control fields,
transferring the computational effort due to such cases to the
preprocessing step: given $\Delta \eps$, the
computational cost of this step only depends 
on the norm $\|\eps\|_{L^\infty(0,T)}$, i.e., on Hypothesis \ref{H}.\\
\end{rem}

\proof
To obtain \eqref{est1} and \eqref{est2}, we will focus first on the local error, i.e. the
approximation obtained on one time step $[t_j,t_{j+1}]$.\\
The sequence $(\psi^K_j)_{j=0,...,N}$ is a time
discretization of the solution of
\begin{eqnarray}\label{psiK}
\left\lbrace
\begin{array}{ll}
i\partial_t \psi^{K}(t) = \left(H_0-\mu\bar\eps(t) \right) \psi^{K} (t),&\bbR^3\times(0,T) \\
\psi^{K}(0)=\psi_0,&\bbR^3
\end{array} \right.
\end{eqnarray}
where the space variable has been omitted and $\bar{\eps}(t) = \bar\eps_{\ell_j}$ is constant over each interval 
$[t_j,t_{j+1}[ = [j\Delta t,(j+1)\Delta t[$, with $j=0,...,N-1$. 
We denote by $\left(S_j(t)\right)_{j=0,...,N-1}$ (instead of
$S_{\ell_j}$) the one-parameter semi-group generated by the operator  
$H_0-\mu\bar\eps_{\ell_j}$ and we  introduce $\delta(t)=\eps(t)-\bar{\eps}$ where $\bar{\eps}$ (instead of $\bar\eps_{\ell_j}$) is the constant value of $\bar{\eps}(t) $ over $[t_j,t_{j+1}]$. 
Therefore, the solution $\psi$ of \eqref{S} is actually the solution of the integral equation, settled for $t\in [t_j,t_{j+1}[$:
\beq\label{id1}
\psi(t) = S_{j}(t-t_j) \psi(t_j) + i\int_{t_j}^{t} S_j(t-s)\mu\delta(s)\psi(s)\ ds.
\eeq
For the upcoming calculations, one should notice that we have 
\begin{equation}
|\delta(t_{j+\frac 12})|\leq\dfrac{\Delta \eps}2\label{estdet}
\end{equation} 
and that for all $t\in [t_j,t_{j+1}] $, 
\beq
|\delta(t)| \leq \dfrac{1}2\left( \Delta \eps+
  \|\dot\eps\|_{L^\infty(0,T)}\Delta t\right).
\label{estde}
\eeq

We consider the following decomposition:
\begin{eqnarray*}
&&\psi(T)-\psi^{K}(T)=\psi(T)-S_{N-1}(\Delta t)\psi(t_{N-1})\\
&&+~\sum_{j=0}^{N-2}S_{N-1}(\Delta t)\dots S_{j+1}(\Delta t) \big(\psi(t_{j+1})
-S_j(\Delta t)\psi(t_j)\big)\\
&&+~S_{N-1}(\Delta t)\dots S_{0}(\Delta t) \psi_0 - \psi^{K}(T)
\end{eqnarray*}
where the last line is equal to $0$ since $\psi^{K}$ satisfies (\ref{psiK}) on $[0,T]$.\\

From now on and in all the following sections, we will consider either that $\|\psi_0\|_{L^2} = 1$ or that $\|\psi_0\|_{H^2} = 1$. From (\ref{SSS}), we know that the operators $S_j$ are isometries in $L^2$. 
Therefore, the use of a triangular inequality brings
\beq\label{eqpsiT}
\|\psi(T)-\psi^{K}(T)\|_{L^2}  \leq  \sum_{j=0}^{N-1} \left\| \psi(t_{j+1})
-S_j(\Delta t)\psi(t_j)\right\|_{L^2}.
\eeq
We will thus calculate and estimate in $L^2$-norm for all $j$ the difference 
\begin{eqnarray}\label{decomp}
\psi(t_{j+1})-S_j(\Delta t)\psi(t_j)&=&  i\int_{t_j}^{t_{j+1}} S_j(t_{j+1}-s)\delta(s)\mu \p(s)\,ds \nonumber \\
&=& i\int_{t_j}^{t_{j+1}} S_j(t_{j+1}-s)\delta(s)\mu \left(\p(s)-S_{j}(s-t_j) \psi(t_j)\right)\,ds \nonumber\\
&&+~ i\int_{t_j}^{t_{j+1}} S_j(t_{j+1}-s)\delta(s)\mu S_{j}(s-t_j) \psi(t_j)\,ds. \nonumber\\
&=& i\int_{t_j}^{t_{j+1}} S_j(t_{j+1}-s)\delta(s)\mu \left(\p(s)-S_{j}(s-t_j) \psi(t_j)\right)\,ds \nonumber\\
&&+~  i\int_{t_j}^{t_{j+1}} \delta(s)\varphi_j(s) \psi(t_j)\,ds
\end{eqnarray} 
where  $\varphi_j(s):=S_j(t_{j+1}-s)\mu S_j(s-t_j)\in \mathcal L(L^2)$.\\

In what follows, we work in parallel on $L^2$ and $H^2$-estimates of $\psi(t_{j+1})-S_j(\Delta t)\psi(t_j)$. We will need basic $L^2$ and $H^2$-estimates of $\p(t)-S_{j}(t-t_j) \psi(t_j)$ for the study of the first integral term of (\ref{decomp}), the second one will be dealt with using a Taylor expansion of $\delta(t)$.\\

From Lemma~\ref{schrod0}, (\ref{id1}) and (\ref{estde}), it is easy to obtain coarse estimates of $\p(t)-S_{j}(t-t_j) \psi(t_j)$. Indeed, for all $t$ in  $[t_j,t_{j+1}] $, one can write

\begin{eqnarray}
 \|\psi(t) - S_{j}(t-t_j) \psi(t_j)\|_{L^2} &=&  
 \left\|i\int_{t_j}^{t} S_j(t-s)\mu\delta(s)\psi(s)\ ds\right\|_{L^2}\nonumber\\
 &\leq&  \int_{t_j}^{t_{j+1}} \left\| S_j(t-s)\mu\delta(s)\psi(s)\right\|_{L^2}\ ds \nonumber\\
 &\leq&\Delta t   \|\mu\|_\nl\dfrac{1}2\left( \Delta \eps +
  \|\dot\eps\|_{L^\infty(0,T)}\Delta t \right) \|\psi_0\|_{L^2}  \nonumber\\
 &\leq&\dfrac{1}2 \|\mu\|_\nl\left( \Delta \eps +
  \|\dot\eps\|_{L^\infty(0,T)}\Delta t \right) \Delta t \label{E1}
  \end{eqnarray} 
and the $H^2$-estimate gives 
\beq
 \|\psi(t) - S_{j}(t-t_j) \psi(t_j)\|_{H^2}  \leq K \|\eps\|_{W^{1,1}(0,T)} \|\mu\|_\nh \left( \Delta \eps +
  \|\dot\eps\|_{L^\infty(0,T)}\Delta t \right) \Delta t \label{E2}
\eeq 
where $K= K(\|\mu\|_{\mathcal L(H^2)}, \eps_{\max})$ is a generic constant that estimates every $\|S_j\|_{\mathcal L(H^2)}$.\\
Therefore, we can obtain more accurate estimates of the first integral
term of \eqref{decomp}.
Thanks to \eqref{E1}, we obtain
\beqn
&&\left\|i\displaystyle\int_{t_j}^{t_{j+1}} S_j(t_{j+1}-s)\delta(s)\mu \left(\p(s)-S_{j}(s-t_j) \psi(t_j)\right)\,ds\right\|_{L^2} \nonumber\\
&&\leq ~\dfrac{1}2\|\mu\|_\nl\left( \Delta \eps +
  \|\dot\eps\|_{L^\infty(0,T)}\Delta t \right) \Delta t \sup_{t\in
  [t_j,t_{j+1}]}\|\p(s)-S_{j}(s-t_j) \psi(t_j)\|_{L^2}\nonumber \\ 
&&\leq~\dfrac{1}4 \|\mu\|^2_\nl  \left( \Delta \eps +
  \|\dot\eps\|_{L^\infty(0,T)}\Delta t \right)^2\Delta t^2 \nonumber \\
&&\leq~\dfrac {1}2\|\mu\|^2_\nl \left( \Delta \eps^2 
+    \|\dot\eps\|^2_{L^\infty(0,T)} \Delta t^2 \right)\Delta t^2.\label{es2}
\eeqn
Working now on the $H^2$-estimate, we deduce from \eqref{E2} in the same way that 
\begin{multline}
\left\|i\displaystyle\int_{t_j}^{t_{j+1}} S_j(t_{j+1}-s)\delta(s)\mu \left(\p(s)-S_{j}(s-t_j) \psi(t_j)\right)\,ds\right\|_{L^2} \\
\leq~\dfrac12 K\|\mu\|_\nh \|\eps\|_{W^{1,1}(0,T)} 
\left(  \Delta \eps^2 +  \|\dot\eps\|_{L^\infty(0,T)}^2\Delta
  t^2 \right)\Delta t^2.\label{es2b} 
\end{multline}
In the two cases ($L^2$ and $H^2$), estimates are stronger than the ones we look
for, and we can focus on the second  integral term of (\ref{decomp}) we want to deal with. \\

We first consider  
\beq\label{phidef}
\begin{array}{llll}
\varphi_j:&[t_j,t_{j+1}]&\to\mathcal L (H^2)\\
&s& \mapsto S_j(t_{j+1}-s)\mu S_j(s-t_j)
\end{array}
\eeq
and note that for all $\psi\in H^2$,~ 
$\|\varphi_j(s) \psi\|_{H^2} = \|S_j(t_{j+1}-s)\mu S_j(s-t_j) \psi\|_{H^2} \leq K\| \psi\|_{H^2}$ 
so that 
\begin{equation}\label{phii}
\forall s\in [t_j,t_{j+1}],\ \|\varphi_j(s)\|_\nh\leq K \|\mu\|_\nh.
\end{equation}
Let us now consider the derivatives of
$\varphi_j(s)$. Since  $(S_j(t))_{t\in \mathbb R}$ denotes the
one-parameter semi-group generated by the operator $H_0-\mu\bar\eps$,
the $\mathcal L(H^2)$ identity $$\partial_t S_j(t) 
=-i\left(H_0-\mu\bar\eps\right)S_j(t)$$ holds and minor calculations give, $\forall s\in [t_j,t_{j+1}]$,
\begin{eqnarray*}
\partial_s \varphi_j(s) &=& iS_j(t_{j+1}-s)[H_0,\mu]S_j(s-t_j)\\
\partial^2_{ss} \varphi_j(s) &=& S_j(t_{j+1}-s) \big[[H_0,\mu], H_0 - \mu\bar\eps \big]S_j(s-t_j).
\end{eqnarray*}
Therefore,
\begin{eqnarray}
\|\partial_s \varphi_j(s)\|_\nh &\leq& K  \|[H_0,\mu]\|_\nh\label{phi}\\
\|\partial^2_{ss} \varphi_j(s)\|_\nh & \leq& K \left\|\big[[H_0,\mu], H_0 - \mu\bar\eps\big]\right\|_\nh.\nonumber
\end{eqnarray}
If we consider the $L^2$-analysis of the method, then $\varphi_j(s) \in \mathcal L(L^2)$ and  $\forall s\in [t_j,t_{j+1}]$, 
\begin{eqnarray}
\|\varphi_j(s)\|_{\mathcal L(L^2)} &\leq&  \|\mu\|_\nl\nonumber\\
\|\partial_s \varphi_j(s)\|_{\mathcal L(L^2)} &\leq&   \|[H_0,\mu]\|_\nl\label{phib}\\
\|\partial^2_{ss} \varphi_j(s)\|_{\mathcal L(L^2)} & \leq&  \left\|\big[[H_0,\mu], H_0 - \mu\bar\eps\big]\right\|_\nl. \nonumber
\end{eqnarray}

Let us now write the third order Taylor expansion of
$t\mapsto \delta(t)=\eps(t)-\bar{\eps}$ in a neighborhood of $t_{j+\frac12}$:
\begin{eqnarray*}
\delta(s)&=&\delta(t_{j+\frac12})+(s-t_{j+\frac12})\dot{\delta}(t_{j+\frac12})+\frac{1}{2}(s-t_{j+\frac12})^2\ddot{\delta}(\theta(s))\\ 
&=&\delta(t_{j+\frac12})+(s-t_{j+\frac12})\dot{\eps}(t_{j+\frac12})+\frac{1}{2}(s-t_{j+\frac12})^2\ddot{\eps}\big(\theta(s)\big),
\end{eqnarray*} 
with $\theta(s)\in [t_j,t_{j+1}]$. We now focus on estimating the term $ i\displaystyle\int_{t_j}^{t_{j+1}} \delta(s)\varphi_j(s) \psi(t_j)\,ds$. 
By means of (\ref{phib}) and
the  $L^2$-norm conservation, we obtain 
$$\left\|\int_{t_j}^{t_{j+1}}  \delta( t_{j+\frac12}) \varphi_j(s)\psi(t_j) ds\right\|_{L^2}
\leq \frac{ 1}{2}\|\mu\|_\nl\Delta \eps \Delta t,$$
\beqn
&&\left\|\int_{t_j}^{t_{j+1}}  \left(s- t_{j+\frac12} \right)\dot{\eps}\left( t_{j+\frac12}\right) \varphi_j(s) \psi(t_j) ds\right\|_{L^2} \nonumber\\
&&= ~\left\|\dot{\eps} \left( t_{j+\frac12} \right)\int_{0}^{
    \frac12{\Delta t}} s \
  \big(\varphi(t_{j+\frac12}+s)-\varphi(t_{j+\frac12}-s)\big) \psi(t_j) ds\right\|_{L^2}\nonumber\\
&&=~\left\|\dot{\eps} \left( t_{j+\frac12}
  \right)\int_{0}^{\frac12{\Delta t}}\int_{t_{j+\frac12}-s}^{t_{j+\frac12}+s}
  s \partial_u \varphi(u) \psi(t_j)\ du ds\right\|_{L^2}\nonumber\\
&&\leq ~\frac1{12}{\|[H_0,\mu]\|_\nl}  \|\dot{\eps}\|_{L^\infty(t_j,t_{j+1})} \Delta t^3\nonumber
\eeqn
and
$$\left\|\int_{t_j}^{t_{j+1}}\frac{1}{2}  \left(s- t_{j+\frac12}
  \right)^2 \ddot{\eps}(\theta(s))\varphi_j(s) \psi(t_j)
  ds\right\|_{L^2} 
\leq\frac1{24}{\|\mu\|_\nl} \|\ddot{\eps}\|_{L^\infty(t_j,t_{j+1})}\Delta t^3.$$
Combining these results with \eqref{es2}, we estimate \eqref{decomp}
  as follows:
\begin{eqnarray*}
\|\psi(t_{j+1})-S_j(\Delta t)\psi(t_j)\|_{L^2} 
&\leq& \frac12 \|\mu\|^2_\nl \left( \Delta \eps^2 \Delta t^2  + \|\dot\eps\|_{\infty}^2 \Delta t^4\right)\\
&& + \frac12\, {\|\mu\|_\nl} \Delta \eps \Delta t\\
&&+~ \frac{1}{24} \left(2\|[H_0,\mu]\|_\nl \|\dot{\eps}\|_{\infty} +\|\mu\|_\nl \|\ddot{\eps}\|_{\infty}\right)\Delta t^3,
\end{eqnarray*}
with $\|\cdot \|_{L^\infty(0,T)} = \| \cdot \|_{\infty}$.
By means of \eqref{eqpsiT}, the global $L^2$-estimate is then:
\begin{eqnarray*}
\|\psi(T)-\psi^{K}(T)\|_{L^2}  
&\leq&  \sum_{j=0}^{N-1} \left\| \psi(t_{j+1})-S_j(\Delta t)\psi(t_j)\right\|_{L^2}\\
&\leq& \frac T2 \|\mu\|^2_\nl \left( \Delta \eps^2 \Delta t  + \|\dot\eps\|_{\infty}^2 \Delta t^3\right) 
+ \frac T2\, {\|\mu\|_\nl} \Delta \eps \\
&&+~ \frac{T}{24} \left(2\|[H_0,\mu]\|_\nl \|\dot{\eps}\|_{\infty} +\|\mu\|_\nl \|\ddot{\eps}\|_{\infty}\right)\Delta t^2,
\end{eqnarray*}
and \eqref{est1} can be deduced with the following constants
$\lambda_1$ and $\lambda_2$ independent of
$\|\eps\|_{L^\infty(0,T)}$~: 
\begin{eqnarray*}
\lambda_1&= &  \frac T2 \|\mu\|_\nl \left( 1+ \Delta \eps \Delta t  \|\mu\|_\nl \right),\\
\lambda_2&=& \dfrac T2\, \|\mu\|_\nl^2 \|\dot \eps\|_{\infty}^2 \Delta t  + \frac T{12}\|[H_0,\mu]\|_\nl
\|\dot{\eps}\|_{\infty} + \frac{T}{24}\|\mu\|_\nl \|\ddot{\eps}\|_{\infty}. 
\end{eqnarray*}

Let us now prove the $H^2$ estimate. 
By means of  \eqref{es2b}, \eqref{phii} and \eqref{phi} and keeping in mind that $K$ is a generic constant 
depending on $\|\mu\|_{\mathcal L(H^2)}$ and  $\eps_{\max}$, we can repeat the previous analysis to find the local estimate:
\begin{eqnarray*}
\|\psi(t_{j+1})-S_j(\Delta t)\psi(t_j)\|_{H^2}  
&\leq&K \|\mu\|_\nh \|\eps\|_{W^{1,1}(0,T)} 
\left(  \Delta \eps^2 \Delta t^2+  \|\dot\eps\|_{L^\infty(0,T)}^2\Delta t^4 \right)\\
 && +~ K  \Delta \eps \Delta t+ K \Big(\|[H_0,\mu]\|_\nh \|\dot{\eps}\|_{L^\infty} + \|\ddot{\eps}\|_{L^\infty}\Big)\Delta t^3.
\end{eqnarray*}
Since one can prove that we can actually write a more precise
estimate of $S_j(\Delta t)$ and replace $K$ by $1+ C \Delta t$ (see properties (\ref{SSS})), we get: 
$$\|S_j(\Delta t)\|_{\mathcal L(H^2)} \leq 1+C \Delta t.$$
and since we have the following intermediate result, where $M>0$
depends on $\|\mu\|_\nh$, $\eps_{\tn{max}}$ and $T$ but is independent
of $N$: 
$$\sum_{j=0}^{N-1} (1+C\Delta t)^{N-j} \Delta t \leq M.$$
The global estimate is obtained as follows:
\begin{eqnarray*}
\|\psi(T)-\psi^{K}(T)\|_{H^2} 
&\leq&  \sum_{j=0}^{N-1} K^{N-j-1}\|\psi(t_{j+1})-S_j(\Delta t)\psi(t_j)\|_{H^2}   \\
&\leq&  \sum_{j=0}^{N-1} (1+C\Delta t)^{N-j} \|\eps\|_{W^{1,1}(0,T)}\left( \Delta \eps^2  \Delta t^2+
  \|\dot\eps\|_{\infty}^2\Delta t^4 \right)\\ 
&&+~  \sum_{j=0}^{N-1} (1+C\Delta t)^{N-j}\left( \Delta \eps \Delta t   + \left( \|[H_0,\mu]\|_\nh \|\dot{\eps}\|_{\infty} +
    \|\ddot{\eps}\|_{\infty} \right)\Delta t^3 \right)\\ 
&\leq&  M \|\eps\|_{W^{1,1}(0,T)}\left( \Delta \eps^2\Delta t +  \|\dot\eps\|_{\infty}^2\Delta t^3 \right) + M \Delta \eps \\
&&+~   M\left( \|[H_0,\mu]\|_\nh
    \|\dot{\eps}\|_{\infty} +  \|\ddot{\eps}\|_{\infty}
  \right)\Delta t^2. 
\end{eqnarray*}
We finally get $\nu_1$ and $\nu_2$ and conclude the proof of Theorem~\ref{anatk}:
\begin{eqnarray*}
\nu_1&= &   M (1 + \|\eps\|_{W^{1,1}(0,T)}T \Delta \eps\Delta t)\\
\nu_2&=&  M \left(\|\eps\|_{W^{1,1}(0,T)} \|\dot\eps\|_{\infty}^2 \Delta t + \|[H_0,\mu]\|_\nh \|\dot{\eps}\|_{\infty} 
+  \|\ddot{\eps}\|_{\infty}\right).
\end{eqnarray*}
\endproof

\begin{rem}
The estimate \eqref{est1} is consistent with the
fact that Algorithm \ref{tk} used with a relevant time discretization is
exact for the piecewise constant control fields.\\
\end{rem}

\section{Improvement in the limit of low intensities}\label{sec:itk}

We now describe a way to improve the time order of the
previous algorithm. Since some constants in the following analysis depend in this
case of the $L^\infty$-norm of the field and the method requires that
the toolkit size scales $\Delta t^3(\eps_{\max}-\eps_{\min})$, it applies in the case of
($L^\infty$-) small control fields. 

\subsection{Algorithm}

The algorithm we propose mixes the toolkit and the splitting approaches, in the
sense that it applies sequentially various operators to correct the
third order local error that appears in the proof of Theorem \ref{anatk}. \\

\begin{algo}(Improved toolkit method for low intensities)\label{Itk}
\begin{enumerate}
\item Preprocessing. Precompute the ``toolkit'', i.e. the set of propagators:
$$S_{\ell}(\Delta t) \textnormal{ for } \ell = 0,\cdots, m,$$
where $(S_{\ell}(t))_{t\in \mathbb R}$ denotes the one-parameter group generated 
by the operator $H_0-\mu\bar{\eps}_\ell$, the sequence $(\eps_\ell)_{\ell=0,\cdots, m,}$ 
being defined by \eqref{tkdisc}.
Include in this set the two special elements:
$$ \Om= e^{\frac{1}{12}[H_0,\mu]\Delta t ^3}, \Theta= e^{\frac{i}{24}\mu\Delta t ^3}$$
and the initial exponents $\alpha_0$ and $\beta_0$ such that ($\eps$ being extended as an even function on [-T,0]): 
\begin{eqnarray*}
\alpha_0 &:=& \frac{\eps(\Delta t)-\eps(0)}{\Delta t} 
= \dot{\eps}(t_{\frac 12}) + \mathcal{O}(\Delta t^2),\\
\beta_0&:=& \frac{\eps(t_{1})-2 \eps(t_{\frac 12})+ \eps(0)}{\Delta t^2}
= \ddot \eps (t_{\frac 12})+\mathcal{O}(\Delta t^2).
\end{eqnarray*}
\item 
Given a control field $\eps\in L^\infty$ satisfying \ref{H} and
  $\psi^{IK}_0=\Om^{\alpha_0}\Theta^{\beta_0}\psi_0$, the sequence $(\psi^{IK}_{j})_{j=0,...,N}$ that
  approximates $(\psi(t_j))_{j=0,...,N}$, is obtained recursively by
  iterating the following loop:
\begin{enumerate} 
\item Find:
$$\ell_j = {\rm argmin}_{\ell=1,\cdots,m}\{| \eps(t_{j+1/2})-\bar{\eps}_\ell|\},$$
\item Compute $\alpha_j$ and $\beta_j$ such that: 
\begin{eqnarray}
\alpha_j &:=& \frac{\eps(t_{j+1})-\eps(t_j)}{\Delta t}\label{dmu}
= \dot{\eps}(t_{j+\frac 12}) + \mathcal{O}(\Delta t^2),\\
\beta_j &:=& \frac{\eps(t_{j+1})-2 \eps(t_{j+\frac 12})+ \eps(t_j)}{\Delta t^2}\label{dnu}
= \ddot \eps (t_{j+\frac 12})+\mathcal{O}(\Delta t^2).
\end{eqnarray}
\item\label{coststep} Set $\psi^{IK}_{j+1}=S_{\ell_j}(\Delta t) \Om^{\alpha_j}\Theta^{\beta_j}\psi^{IK}_j$.\\
\end{enumerate} 
\end{enumerate}
\end{algo}
In many cases, e.g. in the experimental frameworks, only the values of
the field can be handled. The use of exact values for the time
derivatives has then to be avoided when possible. This motivates the
introduction of approximations \eqref{dmu} and
\eqref{dnu} of $\dot{\eps}(t_{j})$
and  $\ddot{\eps}(t_{j})$ in the latest
definitions. The analysis presented hereafter shows that this does not
deteriorate the order of the method.\\

In this method, one must perform two online matrices exponentiations. By
working in a basis where one of these two matrices is diagonal, the
cost of Step~{\it\ref{coststep}} can be reduced to one
exponentiation, making the cost of this method equivalent the second order Strang
operator splitting. 

\subsection{Analysis of the method}

We can now repeat the analysis that has been done in the proof of
Theorem \ref{anatk} to obtain the following estimate.\\

\begin{theorem}\label{anaItk}
Let $\eps\in W^{2,\infty}(0,T)$, $\psi$ be the corresponding solution of
\eqref{S} and $\psi^{IK}$ the approximation of $\psi$ obtained
with Algorithm \ref{Itk}. Given $\Delta t>0$ and $\Delta \eps>0$,
there exists $\lambda'_1>0$, $\lambda'_2>0$, with $\lambda'_1$
independent of $\|\eps\|_{L^\infty(0,T)}$ such that: 
\beq\nonumber
\| \psi(T)-\psi^{IK}(T)\|_{L^2}\leq \lambda'_1 \Delta \eps  + \lambda'_2\Delta t^3.
\eeq
\end{theorem}

\proof
In the framework of this new algorithm, we note that on every time
interval $]t_j,t_{j+1}[$, the approximation $\psi^{IK}$ is the solution of the evolution equation: 
\begin{eqnarray}\label{psiitk}
\left\lbrace
\begin{array}{llll}
i\partial_t \psi^{IK}(t) &= \left(H_0-\mu\bar\eps \right) \psi^{IK} (t),&\bbR^3\times(t_j,t_{j+1}) \\
\psi^{IK}(t_j^+)&=\Om^{\alpha_j}\Theta^{\beta_j}\psi^{IK}(t_j^-) &\bbR^3
\end{array} \right.
\end{eqnarray}
where we set $\psi(0^-) = \psi_0$. We will keep the notations ($S_j$, $\delta (t)$, $\varphi$,...) of the proof of Theorem \ref{anatk}, 
and we first focus on the local error analysis. We consider the following decomposition:
\begin{eqnarray*}
&&\psi(T)-\psi^{IK}(T)=\psi(T)-S_{N-1}(\Delta t)\Om^{\alpha_{N-1}}\Theta^{\beta_{N-1}}\psi(t_{N-1})\\
&&+~\sum_{j=0}^{N-2}S_{N-1}(\Delta t)\Om^{\alpha_{N-1}}\Theta^{\beta_{N-1}}\dots 
S_{j+1}(\Delta t)\Om^{\alpha_{j+1}}\Theta^{\beta_{j+1}} \\
&&\phantom{+~\sum_{j=0}^{N-2}}\times\big(\psi(t_{j+1})
-S_j(\Delta t)\Om^{\alpha_j}\Theta^{\beta_j}\psi(t_j)\big)\\
&&+~S_{N-1}(\Delta t)\Om^{\alpha_{N-1}}\Theta^{\beta_{N-1}}\dots 
S_{0}(\Delta t)\Om^{\alpha_0}\Theta^{\beta_0} \psi_0 - \psi^{IK}(T)
\end{eqnarray*}
where the last line is equal to $0$ since $\psi^{IK}$ satisfies (\ref{psiitk}) on $[0,T]$.\\
The operators $S_j$ are isometries in $L^2$, we will consider that $\|\psi_0\|_{L^2} = 1$ and we also have, for all $j$
\beq\label{omegatheta}
\Om^{\alpha_j}\Theta^{\beta_j} 
=  e^{\frac{\alpha_j}{12}[H_0,\mu]\Delta t ^3} e^{\frac{i\beta_j}{24}\mu\Delta t ^3}  
= \tn{Id} + \left( \frac{\alpha_j}{12}[H_0,\mu] + \frac{i\beta_j}{24}\mu \right) \Delta t ^3   + \tn{Id}~\mathcal O(\Delta t^6)
\eeq
and thus
\beq\label{omegatheta1}
\left\|\Om^{\alpha_j}\Theta^{\beta_j}\right\|_\nl \leq 1+  \mathcal{O}(\Delta t^3).
\eeq
Therefore, the use of a triangular inequality brings
\beq\label{psiT}
\|\psi(T)-\psi^{IK}(T)\|_{L^2}  \leq (1+\mathcal{O}(\Delta t^2))   \sum_{j=0}^{N-1} \left\| \psi(t_{j+1})
-S_j(\Delta t)\Om^{\alpha_j}\Theta^{\beta_j}\psi(t_j)\right\|_{L^2}
\eeq
and we will calculate and estimate in $L^2$-norm for all $j$ the difference 
\begin{eqnarray*}
&&\psi(t_{j+1})-S_j(\Delta t)\Om^{\alpha_j}\Theta^{\beta_j}\psi(t_j)\\
&=& S_j(\Delta t)\psi(t_j) + i\int_{t_j}^{t_{j+1}} S_j(t_{j+1}-s)\delta(s)\mu \p(s)\,ds 
- S_j(\Delta t)\Om^{\alpha_j}\Theta^{\beta_j}\psi(t_j)\\
&=& S_j(\Delta t)\left( \tn{Id} - \Om^{\alpha_j}\Theta^{\beta_j} \right)\psi(t_j) 
+ i\int_{t_j}^{t_{j+1}} S_j(t_{j+1}-s)\delta(s)\mu \p(s)\,ds. \\
\end{eqnarray*} 

We define  $Y(s)=\psi(s)-S_j(s-t_j)\Om^{\alpha_j}\Theta^{\beta_j}\psi(t_j)$ 
for all $s\in[t_j,t_{j+1}]$ and obtain
\begin{multline}
Y(t_{j+1}) = S_j(\Delta t)\left( \tn{Id} - \Om^{\alpha_j}\Theta^{\beta_j} \right)\psi(t_j) \\
+ i\int_{t_j}^{t_{j+1}} S_j(t_{j+1}-s)\delta(s)\mu Y(s)\,ds 
+ i\int_{t_j}^{t_{j+1}} \delta(s)\varphi_j(s)\Om^{\alpha_j}\Theta^{\beta_j}\psi(t_j)\,ds\label{est3}
\end{multline}
where  $\varphi_j(s):=S_j(t_{j+1}-s)\mu S_j(s-t_j)$ and its
derivatives have been estimated in $L^2$ in \eqref{phib}.
As we did in Theorem \ref{anatk}, we start with an estimate of the first integral
term of \eqref{est3}.
For all $t\in]t_j,t_{j+1}]$, we can write:
\begin{eqnarray*}
Y(t)&=&\psi(t)-S_j(t-t_j)\psi(t_j)+S_j(t-t_j)\psi(t_j)-S_j(t-t_j)\Om^{\alpha_j}\Theta^{\beta_j}\psi(t_j)\\
&=&\int_{t_j}^{t}  S_j(t-s)\delta(s)\mu \psi(s)\ ds +S_j(t-t_j)\left(\tn{Id} - \Om^{\alpha_j}\Theta^{\beta_j}\right)\psi(t_j).
\end{eqnarray*}
Moreover, for all $t\in]t_j,t_{j+1}]$, we have
\begin{eqnarray*}
\|Y(t)\|_{L^2}&\leq& \left\|\int_{t_j}^{t}  S_j(t-s)\delta(s)\mu \psi(s)\ ds\right\|_{L^2} + \left\|S_j(t-t_j)\left(\tn{Id} - \Om^{\alpha_j}\Theta^{\beta_j}\right)\psi(t_j)\right\|_{L^2}
\end{eqnarray*}
The operators $S_j$ are isometries in $L^2$ and $\forall t\in [0,T],
\|\p (t)\|_{L^2}=\|\psi_0\|_{L^2}$. Therefore, we deduce from \eqref{omegatheta} that 
$$\left\|S_j(t-t_j)\left(\tn{Id} - \Om^{\alpha_j}\Theta^{\beta_j}\right)\psi(t_j)\right\|_{L^2} \leq 
\left( \frac{\alpha_j}{12}\left\|[H_0,\mu]\right\|_\nl  +\frac{i\beta_j}{24}\left\|\mu\right\|_\nl \right)\Delta t ^3 +  \mathcal{O}(\Delta t^6).$$
Since it is clear that we also have 
$$\left\|\int_{t_j}^{t}  S_j(t-s)\delta(s)\mu \psi(s)\ ds\right\|_{L^2}
\leq   \dfrac 12\|\mu\|_\nl \left( \Delta \eps
+    \|\dot\eps\|_{L^\infty(0,T)} \Delta t \right)\Delta t,$$
one can finally deduce that:
\begin{eqnarray}
&&\left\|i\int_{t_j}^{t_{j+1}}  S_j(t_{j+1}-s)\delta(s)\mu Y(s)\ ds\right\|_{L^2}\nonumber\\
&\leq& \dfrac{1}2\|\mu\|_\nl\left( \Delta \eps +
  \|\dot\eps\|_{L^\infty(0,T)}\Delta t \right) \Delta t \sup_{t\in [t_j,t_{j+1}]}\|Y(t)\|_{L^2} \nonumber\\
 &\leq&~\dfrac 14\|\mu\|^2_\nl \left( \Delta \eps
+    \|\dot\eps\|_{L^\infty(0,T)} \Delta t \right)^2\Delta t^2 +  \mathcal{O}(\Delta \eps \Delta t^4)+  \mathcal{O}(\Delta t^5).\label{supY}
\end{eqnarray}
We focus now on the first and third terms of \eqref{est3}. 
Using \eqref{omegatheta}, we get
$$S_j(\Delta t)\left( \tn{Id} - \Om^{\alpha_j}\Theta^{\beta_j} \right)\psi(t_j) = 
- S_j(\Delta t)\left( \frac{\alpha_j}{12}[H_0,\mu] 
+\frac{i\beta_j}{24}\mu \right)\psi(t_j)\Delta t ^3 + \psi(t_j) \mathcal{O}(\Delta t^6).$$
Let us then consider the second integral term of \eqref{est3}.
On the one hand, we consider the fourth order expansion of $ \delta = \eps - \bar\eps$ 
in a neighborhood of $t_{j+\frac 12}$:
\begin{eqnarray*}
\delta(s)&=&\delta(t_{j+\frac 12})+(s-t_{j+\frac 12})\dot{\delta}(t_{j+\frac 12})+\frac12{(s - t_{j+\frac 12})^2}
\ddot{\delta}(t_{j+\frac 12})+\frac16(s-t_{j+\frac 12})^3\delta^{(3)}(\theta(s))\\  
&=&\delta(t_{j+\frac 12})+(s-t_{j+\frac 12})\dot{\eps}(t_{j+\frac 12})+\frac12{(s-t_{j+\frac 12})^2}
\ddot{\eps}(t_{j+\frac 12})+\frac16{(s-t_{j+\frac 12})^3} \eps ^{(3)}\big(\theta(s)\big)
\end{eqnarray*} 
where $\theta(s)\in [t_j,t_{j+1}]$.
On the other hand, we calculate and/or estimate the four corresponding terms in  
$$ i\int_{t_j}^{t_{j+1}} \delta(s)\varphi_j(s)\Om^{\alpha_j}\Theta^{\beta_j}\psi(t_j)\,ds.$$ 
From \eqref{estdet} and \eqref{phib}, the term of order $0$ gives:
$$\left\|i\int_{t_j}^{t_{j+1}}  \delta(t_{j+\frac 12}) \varphi_j(s) \Om^{\alpha_j}\Theta^{\beta_j}\psi(t_j)ds\right\|_{L^2}
\leq \dfrac 12 \|\mu\|^2_\nl \Delta\eps\Delta t.$$
For the term of order $1$, we can write
\begin{eqnarray*}
&&i\int_{t_j}^{t_{j+1}} (s-t_{j+\frac 12})\dot{\eps}(t_{j+\frac 12}) \varphi_j(s) \Om^{\alpha_j}\Theta^{\beta_j}\psi(t_j) \ ds\\
&=&i~\dot{\eps} \left( t_{j+\frac12} \right)\int_{0}^{\frac12{\Delta t}} s  \big(\varphi_j(t_{j+\frac12}+s)
-\varphi_j(t_{j+\frac12}-s)\big) \Om^{\alpha_j}\Theta^{\beta_j}\psi(t_j) \ ds\\
&=&i~\dot{\eps} \left( t_{j+\frac12} \right)\int_{0}^{\frac12{\Delta t}}\int_{t_{j+\frac12}-s}^{t_{j+\frac12}+s} s \partial_u \varphi_j(u) \Om^{\alpha_j}\Theta^{\beta_j}\psi(t_j)\ du ds\\
&=&i~\dot{\eps}(t_{j+\frac 12})\int_0^{\frac12{\Delta t}} \int_{t_{j+\frac12}-s}^{t_{j+\frac12}+s}s\big(\partial_u\varphi_j(t_{j})
+(u-t_{j})\tau(u) \big) \Om^{\alpha_j}\Theta^{\beta_j}\psi(t_j)\ du ds\\
&=&\frac{\dot{\eps}(t_{j+\frac 12})}{12} S_j\left(\Delta t\right)[H_0,\mu]\Om^{\alpha_j}\Theta^{\beta_j}\psi(t_j)\Delta t^3\\
&& +~i~ \dot{\eps}(t_{j+\frac 12})\int_0^{\frac12{\Delta t}} \int_{t_{j+\frac12}-s}^{t_{j+\frac12}+s}s(u-t_{j}) \tau(u-t_j)\Om^{\alpha_j}\Theta^{\beta_j}\psi(t_j)\ du ds\\
&=&\frac{ \alpha_j}{12} S_j\left(\Delta t\right)[H_0,\mu]\psi(t_j)\Delta t^3 
+ S_j\left(\Delta t\right)[H_0,\mu]\psi(t_j) \mathcal O(\Delta t^6)\\
&& +~i~\dot{\eps}(t_{j+\frac 12}) \int_0^{\frac12{\Delta t}} \int_{t_{j+\frac12}-s}^{t_{j+\frac12}+s}s(u-t_{j}) \tau(u-t_j)\Om^{\alpha_j}\Theta^{\beta_j}\psi(t_j)\ du ds
\end{eqnarray*} 
where we used \eqref{dmu}, \eqref{omegatheta} and \eqref{phib} and the function $\tau : s\in[0,\Delta t] \mapsto \tau(s)\in \mathcal
L(L^2)$ is defined as the function that appears in the following
expansion of $\partial_u \varphi_j$ around $t_{j}$, for any $\psi\in
L^2$ 
\begin{eqnarray*}
\partial_u\varphi_j(u)\psi&=& \partial_u\varphi_j(t_j)\psi+(u-t_j)\tau(u-t_j)\psi\\
&=&iS_j( \Delta t)[H_0,\mu]\psi+(u-t_j) \tau(u-t_j)\psi.
\end{eqnarray*} 
Using the estimate (coming from \eqref{phib})
\beq\label{depc}
\|\tau (s)\|_\nl\leq \left\| \Big[[H_0,\mu],H_0-\mu\bar\eps \Big]\right\|_\nl ~~~\forall s\in[0,\Delta t],
\eeq
along with $\|\p(t_j)\|_{L^2} =  \|\p_0\|_{L^2} = 1$, \eqref{dmu} and \eqref{omegatheta1} we find that for all $j$,
\begin{eqnarray*}
&&\left\|i \dot{\eps}(t_{j+\frac 12})\int_0^{\frac12{\Delta t}} \int_{t_{j+\frac12}-s}^{t_{j+\frac12}+s}s(u-t_{j}) \tau(u-t_j)\Om^{\alpha_j}\Theta^{\beta_j}\psi(t_j)\ du ds\right\|_{L^2}\\
&\leq& \left( \alpha_j +\mathcal{O}\left(\Delta t^2\right) \right)
\int_0^{\frac12{\Delta t}} \int_{t_{j+\frac12}-s}^{t_{j+\frac12}+s}s(u-t_{j}) \left\|\Big[[H_0,\mu],H_0-\mu\bar\eps \Big]\right\|_\nl\left\|\psi(t_j)\right\|_{L^2}\ du ds. \\
&\leq& \dfrac{\alpha_j}{24}\left\|\Big[[H_0,\mu],H_0-\mu\bar\eps \Big]\right\|_\nl\Delta t^4    + \mathcal{O}\left(\Delta t^5\right). 
\end{eqnarray*}
We also prove easily that for all $j$, 
$$\left\|S_j\left(\Delta t\right)[H_0,\mu]\psi(t_j) \mathcal O(\Delta t^6)\right\|_{L^2}= \mathcal{O}\left(\Delta t^6\right).$$
For the term of order $2$, using \eqref{dnu}, \eqref{omegatheta} and \eqref{phib} and the first order expansion of $\varphi_j$ around $t_{j}$, $\varphi_j(s)\psi = \varphi_j(t_j)\psi +(s-t_j)\theta(s-t_j)\psi$ for all $\psi\in L^2$, defining $\theta : s\in[0,\Delta t] \mapsto \theta(s)\in \mathcal L(L^2)$, we can write
\begin{eqnarray*}
&&i\int_{t_j}^{t_{j+1}}\frac12{(s-t_{j+\frac 12})^2}\ddot{\eps}(t_{j+\frac 12})\varphi_j(s)\Om^{\alpha_j}\Theta^{\beta_j}\psi(t_j)\ ds\\
&=&i\ddot{\eps}(t_{j+\frac 12})~\int_{t_j}^{t_{j+1}}\frac12{(s-t_{j+\frac 12})^2}\left( \varphi_j(t_j) +(s-t_j)\theta(s-t_j)  \right) \Om^{\alpha_j}\Theta^{\beta_j}\psi(t_j)\ ds\\
&=&\frac{i\ddot{\eps}(t_{j+\frac 12})}{24} S_j\left(\Delta t\right)\mu\Om^{\alpha_j}\Theta^{\beta_j}\psi(t_j)\Delta t^3\\
&& +~\frac{i \ddot{\eps}(t_{j+\frac 12})}2 \int_{t_j}^{t_{j+1}}(s-t_{j+\frac 12})^2(s-t_j)\theta(s-t_j)\Om^{\alpha_j}\Theta^{\beta_j}\psi(t_j)\ ds\\
&=&\frac{i\beta_j }{24} S_j\left(\Delta t\right)\mu\psi(t_j)\Delta t^3 
+ S_j\left(\Delta t\right)\mu \psi(t_j) \mathcal O(\Delta t^6)\\
&& +~\frac{i \ddot{\eps}(t_{j+\frac 12})}2 \int_{t_j}^{t_{j+1}}(s-t_{j+\frac 12})^2(s-t_j)\theta(s-t_j)\Om^{\alpha_j}\Theta^{\beta_j}\psi(t_j)\ ds.
\end{eqnarray*} 
Using \eqref{phib}, we get the estimate $\|\theta (s)\|_\nl\leq \left\| \Big[H_0,\mu]\right\|_\nl $, $\forall s\in[0,\Delta t],$ and using it with \eqref{dnu} and \eqref{omegatheta1}, we obtain that for all $j$,
\begin{eqnarray*}
&&\left\|\frac{i \ddot{\eps}(t_{j+\frac 12})}2 \int_{t_j}^{t_{j+1}}(s-t_{j+\frac 12})^2(s-t_j)\theta(s-t_j)
\Om^{\alpha_j}\Theta^{\beta_j}\psi(t_j)\ ds\right\|_{L^2}\\
&\leq&\frac 12 \left( \beta_j +\mathcal{O}\left(\Delta t^2\right) \right)
\int_{t_j}^{t_{j+1}}(s-t_{j+\frac 12})^2(t_j-s) \left\|[H_0,\mu]\right\|_\nl\left\|\psi(t_j)\right\|_{L^2}\  ds \\
&\leq& \frac 12 \left( \beta_j +\mathcal{O}\left(\Delta t^2\right) \right)\left\|[H_0,\mu]\right\|_\nl
\int_{-\frac{\Delta t}2}^{\frac{\Delta t}2}u^2\left(\frac{\Delta t}2 -u\right) \  du\\
&\leq&\dfrac{\beta_j}{48} \left\|[H_0,\mu]\right\|_\nl \Delta t^4    + \mathcal{O}\left(\Delta t^5\right). 
\end{eqnarray*}
and we also prove easily that $\left\|S_j\left(\Delta t\right)\mu\psi(t_j) \mathcal O(\Delta t^6)\right\|_{L^2}= \mathcal{O}\left(\Delta t^6\right).$\\

Combining these results with \eqref{supY} into equation \eqref{est3}, we obtain:
\begin{eqnarray*}
\|Y(t_{j+1}) \|_{L^2} &=&  \left\| \psi(t_{j+1})
-S_j(\Delta t)\Om^{\alpha_j}\Theta^{\beta_j}\psi(t_j)\right\|_{L^2}\\
&\leq&
\left\|S_j(\Delta t)\left( \tn{Id} - \Om^{\alpha_j}\Theta^{\beta_j} \right)\psi(t_j) 
+ i\int_{t_j}^{t_{j+1}} \delta(s)\varphi_j(s)\Om^{\alpha_j}\Theta^{\beta_j}\psi(t_j)\,ds\right\|_{L^2}\\
&&+~ \left\|i\int_{t_j}^{t_{j+1}} S_j(t_{j+1}-s)\delta(s)\mu Y(s)\,ds \right\|_{L^2}\\
&\leq& \dfrac 12 \|\mu\|^2_\nl \Delta\eps\Delta t + \dfrac{\alpha_j}{24}\left\|\Big[[H_0,\mu],H_0-\mu\bar\eps \Big]\right\|_\nl\Delta t^4   \\
&&+~\dfrac{\beta_j}{48} \left\|[H_0,\mu]\right\|_\nl \Delta t^4  +\dfrac 12\|\mu\|^2_\nl  \left(\Delta \eps^2\Delta t^2+ \|\dot\eps\|^2_{L^\infty(0,T)} \Delta t^4\right)   \\
&&+~ \mathcal{O}(\Delta \eps \Delta t^4)+  \mathcal{O}(\Delta t^5)
\end{eqnarray*} 
We have now a local in time estimate that should be traduced in a global one, and from \eqref{psiT}, we get
\begin{eqnarray*}
&&\|\psi(T)-\psi^{IK}(T)\|_{L^2} \\
& \leq &(1+\mathcal O(\Delta t^2)) \sum_{j=0}^{N-1}  \left\| \psi(t_{j+1})
-S_j(\Delta t)\Om^{\alpha_j}\Theta^{\beta_j}\psi(t_j)\right\|_{L^2}\\
& \leq & \dfrac T2 \|\mu\|^2_\nl \Delta\eps + \dfrac{\alpha_j T}{24}\left\|\Big[[H_0,\mu],H_0-\mu\bar\eps \Big]\right\|_\nl\Delta t^3   \\
&&+~\dfrac{\beta_j T}{48} \left\|[H_0,\mu]\right\|_\nl \Delta t^3  +\dfrac T2\|\mu\|^2_\nl  \left(\Delta \eps^2\Delta t+ \|\dot\eps\|^2_{L^\infty(0,T)} \Delta t^3\right)   \\
&&+~ \mathcal{O}(\Delta \eps \Delta t^3)+  \mathcal{O}(\Delta t^4).
\end{eqnarray*}
The result follows, with 
$$\lambda'_1  = \frac T2  \|\mu\|^2_\nl(1+\Delta \eps \Delta t ) $$ and  
$$\lambda'_2 = \dfrac{\alpha_j T}{24}\left\|\Big[[H_0,\mu],H_0-\mu\bar\eps \Big]\right\|_\nl
+\dfrac{\beta_j T}{48} \left\|[H_0,\mu]\right\|_\nl   +\dfrac T2\|\mu\|^2_\nl \|\dot\eps\|^2_{L^\infty(0,T)}$$  
\endproof


In this theorem, the constants $\lambda'_2$ depends on $\| \eps\|_{L^\infty(0,T)}$
through the commutator $\Big[[H_0,\mu],H_0-\mu\bar\eps \Big]$ that appears in
\eqref{depc}. This contrasts with the result obtained in Theorem~\ref{anatk}. The
explanation of this situation comes from the fact that the norms of
$\varphi_j(s):=S_j(t_{j+1}-s)\mu S_j(s-t_j)$ (defined in \eqref{phidef}) 
and its first derivative does not depend on $\| \eps\|_{L^\infty(0,T)}$, whereas its
second derivative does. Thus, errors in Algorithm~\ref{Itk} 
depend on $L^\infty$-norm of the control field as in the case of the second order
Strang operator splitting. Although these two methods present the same
computational complexity, the order of Algorithm~\ref{Itk} is
higher when $\Delta \eps$ scales $\Delta t^3$.
~

\section{Improvement in the limit of large intensities}\label{sec:itk2}

We now describe a way to improve the time order of the
Algorithm \ref{tk} in the case of large intensities. The following method
 enables to replace $\Delta \eps$ by $\Delta \eps\Delta
 t $ in the estimates. 

\subsection{Algorithm}

The algorithm we propose improve the accuracy in the approximation of
$\varepsilon$. This improvement is obtained by using two toolkit
elements instead of one at each time step. \\

\begin{algo}(Improved toolkit method for large intensities)\label{Itk2}
\begin{enumerate}
\item Preprocessing. Precompute the ``toolkit'', i.e. the set of propagators:
$$S_{\ell}(\Delta t) \textnormal{ for } \ell = 0,\cdots, m,$$
where $(S_{\ell}(t))_{t\in \mathbb R}$ denotes the one-parameter group
generated by the operator $H_0-\mu\bar{\eps}_\ell$, the sequence
$(\eps_\ell)_{\ell=0,\cdots, m,}$ being defined by \eqref{tkdisc}. 
\item 
Given a control field $\eps\in L^\infty$ satisfying \ref{H} and
  $\psi^{JK}_0=\psi_0$, the sequence $(\psi^{JK}_{j})_{j=0,...,N}$ that
  approximates $(\psi(t_j))_{j=0,...,N}$, is obtained recursively by
  iterating the following loop:
\begin{enumerate} 
\item Find $\ell_j$ such that:
$$\eps(t_{j+1/2})\in [\bar{\eps}_{\ell_j},\bar{\eps}_{\ell_j+1}].$$
\item Compute $\alpha_j$ and $\beta_j$ such that: 
\begin{eqnarray}
\alpha_j\bar{\eps}_{\ell_j}+\beta_j\bar{\eps}_{\ell_j+1} &=&\eps(t_{j+1/2})\nonumber\\
\alpha_j+\beta_j &=& 1 \label{alphabeta2}
\end{eqnarray}
\item\label{coststep2} Set $\psi^{JK}_{j+1}=S_{\ell_j+1}(\Delta t )^{\beta_j}S_{\ell_j}(\Delta t)^{\alpha_j}\psi^{JK}_j$.\\
\end{enumerate} 
\end{enumerate}
\end{algo}
In this method, one must perform two online matrices
exponentiations. The cost of the corresponding step, namely Step~{\it\ref{coststep2}} can be reduced to three matrix
products when precomputing the mappings between the diagonalization basis of two
consecutive toolkit elements.
\begin{rem}\label{quantif}
Another way to reduce the cost of this step, is to quantify the
values of $\alpha_j$ (and $\beta_j$) and precompute a toolkit
containing elements of the form : $S_{\ell_j+1}(\Delta t
)^{\beta_j}S_{\ell_j}(\Delta t)^{\alpha_j}$. This method is tested
in Sec. \ref{sec:nr}.
\end{rem}
\subsection{Analysis of the method}

We can now repeat the analysis that has been done in the proof of
Theorem \ref{anatk} to obtain the following estimate.\\

\begin{theorem}\label{anaItk2}
Let $\eps\in W^{3,\infty}(0,T)$, $\psi$ be the corresponding solution of
\eqref{S} and $\psi^{JK}$ the approximation of $\psi$ obtained
with Algorithm \ref{Itk}. Given $\Delta t>0$ and $\Delta \eps>0$,
there exists $\lambda''_1>0$, $\lambda''_2>0$, both independent of $\|\eps\|_{L^\infty(0,T)}$ such that: 
\beq\nonumber
\| \psi(T)-\psi^{JK}(T)\|_{L^2}\leq  \lambda''_1 \Delta \eps\Delta t + \lambda''_2\Delta t^2 .
\eeq
\end{theorem}

\proof
In this algorithm, two control fields are involved successively in the propagation over the
interval $[t_j,t_{j+1}]$. As in the previous proofs, we introduce
$\delta(s)=\varepsilon(s)-\bar{\varepsilon}(s)$, with
 $$\bar{\varepsilon}(s)=\left\{ \begin{array}{cl}\bar{\varepsilon}_{\ell_j}   & s\in [t_j,t_j+\alpha_j\Delta t[,\\
                                                 \bar{\varepsilon}_{\ell_j+1} & s\in [t_j+\alpha_j\Delta t,t_{j+1}[. 
                                                  \end{array}\right.$$
Note first that for all $s\in [t_j,t_{j+1}]$ 
\beq\label{estimdelta}
|\delta(s)|\leq \Delta \eps + \frac12\|\dot{\varepsilon}\|_{L^\infty(0,T)} \Delta t
\eeq
and denote by $\left(S_j(t)\right)_{j=0,...,N-1}$ and $\left(S'_j(t)\right)_{j=0,...,N-1}$ 
the one-parameter semi-groups generated by the operators $H_0-\mu\bar\eps_{\ell_j}$ 
and $H_0-\mu\bar\eps_{\ell_j+1}$  respectively.\\
 Following the same analysis as for Algorithm~\ref{tk}, we set 
 $(\psi^{JK}_j)_{j=0,...,N}$ as the time discretization of the solution of:
\begin{eqnarray}\label{psiJK}
\left\lbrace
\begin{array}{ll}
i\partial_t \psi^{JK}(t) = \left(H_0-\mu\bar\eps(t) \right) \psi^{JK} (t),&\bbR^3\times(0,T) \\
\psi^{JK}(0)=\psi_0,&\bbR^3
\end{array} \right.
\end{eqnarray}
where $\bar{\eps}(t)$ (defined right above) is constant over each interval 
$[t_j,t_j+\alpha_j\Delta t[ $ and $ [t_j+\alpha_j\Delta t,t_{j+1}[$, with $j=0,...,N-1$. 
In the same way as we obtained \eqref{id1}, the solution $\psi$ of \eqref{S} satisfies,
$$\p(t_j+\alpha_j\Delta t) = S_{j}(\alpha_j \Delta t) \p(t_j)
-i\int_{t_j}^{t_{j}+\alpha_j\Delta t} S_j(t_j+\alpha_j\Delta t - s)\mu\delta(s)\psi(s)\ ds$$
and
$$\p(t_{j+1}) = S'_{j}(\beta_j \Delta t) \p(t_j+\alpha_j \Delta t)
-i\int_{t_j+\alpha_j\Delta t}^{t_{j+1}} S'_j(t_{j+1} - s)\mu\delta(s)\psi(s)\ ds.$$
As in (\ref{decomp}), it gives rise to:
\begin{eqnarray}
&&\p(t_{j+1})- S'_j(\beta_j\Delta t)S_j(\alpha_j\Delta t)\p(t_j) \nonumber\\
&=&i\int_{t_j+\alpha_j\Delta t}^{t_{j+1}} S'_j(t_{j+1}-s)\mu\delta(s) \p(s) \ ds \nonumber\\ 
&&+  i\int_{t_j}^{t_{j}+\alpha_j\Delta t} S'_j(\beta_j\Delta t)S_j(t_{j}+\alpha_j\Delta t-s)\mu\delta(s) \p(s)\ ds \nonumber\\ 
&=&i\int_{t_j+\alpha_j\Delta t}^{t_{j+1}} S'_j(t_{j+1}-s)\mu\delta(s) \left(\p(s) - S'_j(s-t_j-\alpha_j\Delta t) \p(t_j+\alpha_j\Delta t)\right)    \ ds \nonumber\\ 
&&+  i\int_{t_j}^{t_{j}+\alpha_j\Delta t} S'_j(\beta_j\Delta t)S_j(t_{j}+\alpha_j\Delta t-s)\mu\delta(s)\left(\p(s) - S_j(s-t_j) \p(t_j)\right)\ ds \nonumber\\ 
&& +i\int_{t_j+\alpha_j\Delta t}^{t_{j+1}}   \delta(s)\widetilde{\varphi}'_j(s)S_j(\alpha_j\Delta t) \psi(t_j) \ ds \nonumber\\
 &&+ i\int_{t_j}^{t_{j}+\alpha_j\Delta t} S'_j(\beta_j \Delta   t)    \delta(s)\widetilde{\varphi}_j(s) \psi(t_j) \ ds .\label{id5}
\end{eqnarray} 
where
$\widetilde{\varphi}'_j(s):=S'_j(t_{j+1}-s)\mu S'_j(s-t_j-\alpha_j\Delta t)$ and $\widetilde{\varphi}_j(s):=S_j(t_j+\alpha_j\Delta t-s)\mu S_j(s-t_j)$.
As in the proof of Theorem~\ref{anatk} (see right above \eqref{eqpsiT}) we use the appropriate decomposition
\begin{eqnarray*}
&&\psi(T)-\psi^{JK}(T)=\psi(T)-S'_{N-1}(\beta_{N-1}\Delta t)S_{N-1}(\alpha_{N-1}\Delta t)\psi(t_{N-1})\\
&&+~\sum_{j=0}^{N-2}S'_{N-1}(\beta_{N-1}\Delta t)S_{N-1}(\alpha_{N-1}\Delta t)\dots S'_{j+1}(\beta_{j+1}\Delta t)S_{j+1}(\alpha_{j+1}\Delta t)\\
&&\hspace{7cm} \times \big(\psi(t_{j+1})-S'_{j}(\beta_j\Delta t)S_{j}(\alpha_j\Delta t)\psi(t_j)\big)\\
&&+~S'_{N-1}(\beta_{N-1}\Delta t)S_{N-1}(\alpha_{N-1}\Delta t)\dots S'_0(\beta_0\Delta t)S_0(\alpha_0\Delta t) \psi_0 - \psi^{JK}(T)
\end{eqnarray*}
where the last line is equal to $0$ since $\psi^{JK}$ satisfies (\ref{psiJK}) on $[0,T]$.
We have the corresponding estimate (see \eqref{eqpsiT})
$$\|\psi(T)-\psi^{JK}(T)\|_{L^2}  \leq  \sum_{j=0}^{N-1} \left\| \p(t_{j+1})- S'_j(\beta_j\Delta t)S_j(\alpha_j\Delta t)\p(t_j)\right\|_{L^2}$$
and we will thus calculate and estimate in $L^2$-norm for all $j$ the four terms of  \eqref{id5}.
As in \eqref{es2}, but using now the new estimate \eqref{estimdelta} of $\delta$, the two first terms of the right hand side of \eqref{id5} can be respectively estimated by:
\begin{multline}
\left\|i\int_{t_j+\alpha_j\Delta t}^{t_{j+1}} S'_j(t_{j+1}-s)\mu\delta(s) \left(\p(s) 
- S'_j(s-t_j-\alpha_j\Delta t) \p(t_j+\alpha_j\Delta t)\right)    \ ds\right\|_{L^2}\\
\leq~\beta_j\|\mu\|^2_\nl \left( \Delta \eps^2 \Delta t^2 + \frac12\|\dot\eps\|^2_{L^\infty(0,T)} \Delta t^4 \right)\label{t0}
\end{multline}
and
\begin{multline}
\left\|  i\int_{t_j}^{t_{j}+\alpha_j\Delta t} S'_j(\beta_j\Delta t)S_j(t_{j}+\alpha_j\Delta t-s)\mu\delta(s)\left(\p(s) - S_j(s-t_j) \p(t_j)\right) \right\|_{L^2}\\
\leq~\alpha_j\|\mu\|^2_\nl\left( \Delta \eps^2 \Delta t^2 + \frac12\|\dot\eps\|^2_{L^\infty(0,T)} \Delta t^4 \right).\label{t1}
\end{multline}
Let us now focus on the third and fourth terms of \eqref{id5}. We have: 
\beqn
\int_{t_j+\alpha_j\Delta t}^{t_{j+1}}
      \delta(s)\widetilde{\varphi}'_j(s)S_j(\alpha_j\Delta t)
      \psi(t_j)     \ ds 
=\int_{t_j+\alpha_j\Delta t}^{t_{j+1}}
      \delta(s)\widetilde{\varphi}'_j(t_j+\alpha_j\Delta t)S_j(\alpha_j\Delta t)
      \psi(t_j)     \ ds \nonumber\\
+\int_{t_j+\alpha_j\Delta t}^{t_{j+1}}
      \delta(s)\int_{t_j+\alpha_j\Delta t}^s\partial_u\widetilde{\varphi}'_j(u)
          \ du\ S_j(\alpha_j\Delta t)\psi(t_j) \ ds\nonumber\\
=\int_{t_j+\alpha_j\Delta t}^{t_{j+1}}
      \delta(s)S'_j(\beta_j \Delta t)\mu S_j(\alpha_j\Delta t)
      \psi(t_j)     \ ds \nonumber\\
+\int_{t_j+\alpha_j\Delta t}^{t_{j+1}}
      \delta(s)\int_{t_j+\alpha_j\Delta
        t}^s\partial_u\widetilde{\varphi}'_j(u)\ du  \ S_j(\alpha_j\Delta t)
      \psi(t_j)      \ ds\nonumber
\eeqn
and
\beqn
\int_{t_j}^{t_{j}+\alpha_j\Delta t} S'_j(\beta_j \Delta
      t)    \delta(s)\widetilde{\varphi}_j(s) \psi(t_j) \ ds=\int_{t_j}^{t_{j}+\alpha_j\Delta t}\delta(s) S'_j(\beta_j \Delta
      t)    \widetilde{\varphi}_j(t_{j}+\alpha_j\Delta t)
      \psi(t_j) \ ds \nonumber\\
-\int_{t_j}^{t_{j}+\alpha_j\Delta t} \delta(s) S'_j(\beta_j \Delta
      t)   \int_{s}^{t_j+\alpha_j\Delta
        t}\partial_u\widetilde{\varphi}_j(u) \ du\ \psi(t_j) \ ds \nonumber\\
=\int_{t_j}^{t_{j}+\alpha_j\Delta t} \delta(s) S'_j(\beta_j \Delta
      t)   \mu S_j(\alpha_j\Delta t)
      \psi(t_j) \ ds \nonumber\\
-\int_{t_j}^{t_{j}+\alpha_j\Delta t}    \delta(s) S'_j(\beta_j \Delta
      t)\int_{s}^{t_j+\alpha_j\Delta
        t} \partial_u\widetilde{\varphi}_j(u) \ du \ \psi(t_j)\ ds \nonumber.
\eeqn
By means of \eqref{alphabeta2}, we have:
\beqn
\int_{t_j}^{t_{j+1}}\delta (s)\ ds& =&\int_{t_j}^{t_{j+1}} \varepsilon(s)-\varepsilon(t_{j+1/2})\ ds
+ \int_{t_j}^{t_{j+1}}\varepsilon(t_{j+1/2})- \bar{\varepsilon}(s)\
ds\nonumber\\
& =&\int_{t_j}^{t_{j+1}} \varepsilon(s)-\varepsilon(t_{j+1/2})\ ds\nonumber\\
&= &\int_{t_j}^{t_{j+1}}\ddot{\varepsilon}(\theta(s))\frac12(s-t_{j+1/2})^2\ ds,\nonumber
\eeqn
where $\theta(s)\in [t_j,t_{j+1}]$. Consequently,
\beq
\left\|\int_{t_j}^{t_{j+1}}
      \delta(s)S'_j(\beta_j \Delta t)\mu S_j(\alpha_j\Delta t)
      \psi(t_j)     \ ds \right\|_{L^2}\leq \frac1{24}\|\mu\|_\nl\|\ddot{\varepsilon}\|_{L^\infty(0,T)}\Delta t^3.\label{t2}
\eeq
From \eqref{phib} and \eqref{estimdelta}, we obtain
\begin{multline}
\left\| \int_{t_j+\alpha_j\Delta t}^{t_{j+1}} \delta(s)\int_{t_j+\alpha_j\Delta  t}^s\partial_u\widetilde{\varphi}'_j(u)\ du  \ S_j(\alpha_j\Delta t)   \psi(t_j)\ ds\right\|_{L^2}\\
\leq \frac12\beta_j^2\|[H_0,\mu]\|_\nl \left(\Delta
  \varepsilon + \frac12 \|\dot{\varepsilon}\|_{L^\infty(0,T)}\Delta t  \right)\Delta t^2\label{t3}
\end{multline}
and similarly, we find that:
\begin{multline}
\left\| \int^{t_j+\alpha_j\Delta t}_{t_{j}}
      \delta(s)S'_j(\beta_j\Delta t)\int^{t_j+\alpha_j\Delta
        t}_s\partial_u\widetilde{\varphi}_j(u)\ du  \
      \psi(t_j)\ ds   \right\|_{L^2}\\
      \leq \frac12\alpha_j^2\|[H_0,\mu]\|_\nl \left( \Delta
  \varepsilon + \frac12 \|\dot{\varepsilon}\|_{L^\infty(0,T)}\Delta t \right)\Delta t^2.\label{t4}
\end{multline}
Combining \eqref{t0}, \eqref{t1}, \eqref{t2}, \eqref{t3} and \eqref{t4}, we obtain:
\begin{eqnarray*}
\left\|\p(t_{j+1})- S'_j(\beta_j\Delta t)S_j(\alpha_j\Delta t) \p(t_j)\right\|_{L^2}
&\leq& \|\mu\|^2_\nl \left( \Delta \eps^2+ \frac12\|\dot\eps\|^2_{L^\infty(0,T)} \Delta t^2\right) \Delta t^2\\
&&+ \frac1{24}\|\mu\|_\nl\|\ddot{\varepsilon}\|_{L^\infty(0,T)}\Delta t^3\\
&&+\frac12\|[H_0,\mu]\|_\nl \left( \Delta \varepsilon + \frac12 \|\dot{\varepsilon}\|_{L^\infty(0,T)}\Delta t \right)\Delta t^2.
\end{eqnarray*}
The global estimate follows
\begin{eqnarray*}
\|\psi(T)-\psi^{JK}(T)\|_{L^2}  &\leq&  \sum_{j=0}^{N-1} \left\| \p(t_{j+1})- S'_j(\beta_j\Delta t)S_j(\alpha_j\Delta t)\p(t_j)\right\|_{L^2}\\
&\leq & \|\mu\|^2_\nl \left( \Delta \eps^2+ \frac12\|\dot\eps\|^2_{L^\infty(0,T)} \Delta t^2\right)T \Delta t\\
&&+ \frac1{24}\|\mu\|_\nl\|\ddot{\varepsilon}\|_{L^\infty(0,T)}T\Delta t^2\\
&&+\frac12\|[H_0,\mu]\|_\nl \left( \Delta \varepsilon + \frac12 \|\dot{\varepsilon}\|_{L^\infty(0,T)}\Delta t \right)T\Delta t
\end{eqnarray*}
and the proof of Theorem~\ref{anaItk2} is complete, with
$$\lambda''_1=\frac12\|[H_0,\mu]\|_\nl T + \|\mu\|_\nl^2T\Delta \eps,$$
\begin{eqnarray*}
\lambda''_2&=&\frac14\|[H_0,\mu]\|_\nl\|\dot{\varepsilon}\|_{L^\infty(0,T)}T+\frac1{24}\|\mu\|_\nl\|\ddot{\varepsilon}\|_{L^\infty(0,T)}T
\\&&+\frac12\|\mu\|_\nl^2\|\dot{\varepsilon}\|_{L^\infty(0,T)}T\Delta
t.
\end{eqnarray*}

\endproof
\section{Numerical results}\label{sec:nr}
In this section, we check numerically that the order of the estimates we
have obtained in this paper are optimal, and we compare computational
costs of the methods.

\subsection{Model}

In order to test the performance of the algorithms on
a realistic case, a model already
treated in the literature has been considered. The system is a molecule of HCN
modeled as a rigid rotator. We refer the reader to \cite{claudedionjcp} for numerical details
concerning this system.\\
As a control field, we use an arbitrary field of the form
$\varepsilon(t)=\varepsilon_{\max}\sin(\omega t)$, with
$\varepsilon_{\max}=5.10^{-5}$ and $\omega= 5.10^{-6}$. The
parameters are chosen in accordance with usual scales considered for
this model. The use of an analytic
formula for the field enables us to work with exact values,
i.e. to test the cases $\Delta \varepsilon=0$.

\subsection{Orders of convergence}

\begin{centering}
\begin{figure}[h!]
\psfrag{dt}[r][c]{ $\Delta t/T$}
\psfrag{error}[c][t]{ $\| \psi(T)-\psi^{num}(T)\|_{L^2}$}
\center
\includegraphics[width=.8\textwidth]{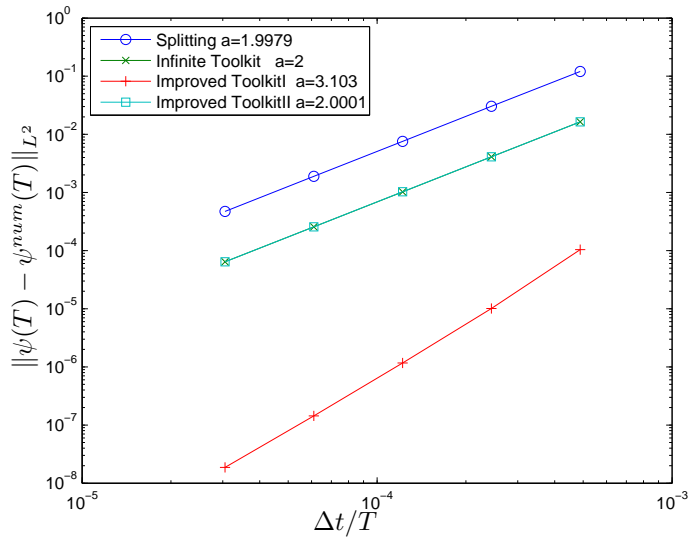}
\caption{Error with respect to $\Delta t$, when $\Delta
  \varepsilon=0$ for toolkit method and Improved toolkit I method, and
  when$\Delta
  \varepsilon=c\Delta t$ for Improved toolkit II method. Here, $\psi^{num}$ stands for the approximation
  of $\psi$ when using the toolkit method, the second order Strang
operator splitting, the Improved toolkit I method and the Improved
toolkit II.  The coefficient $a$ is the regression
coefficient.}\label{time_error} 
\end{figure}
\end{centering}

To test the time order, we first work with $\Delta \varepsilon=0$,
with various values of $\Delta t$. The numerical orders correspond to
the ones obtained in our analysis. Curves of convergence are depicted
in Fig. \ref{time_error}. 

\begin{centering}
\begin{figure}[h!]
\psfrag{dc}[c][b]{ $\Delta \varepsilon/\varepsilon_{\max}$}
\psfrag{error}[c][t]{ $\| \psi(T)-\psi^{num}(T)\|_{L^2}$}
\center
\includegraphics[width=.8\textwidth]{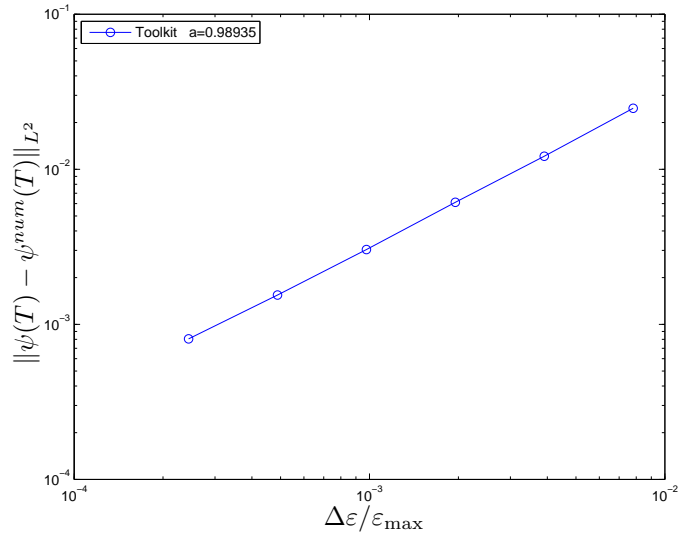}
\caption{Error with respect to $\Delta \varepsilon$, when $\Delta t$  is small. Here, $\psi^{num}$ stands for the approximation  of $\psi$ when using the toolkit method.   }\label{toolkit_error}
\end{figure}
\end{centering}

The order with respect to $\Delta \varepsilon$ is also obtained
numerically by using a small time step. In this test, the numerical
order is consistent 
with the one obtained in Theorem~\ref{anatk}. The convergence with
respect to this parameter is presented in Fig. \ref{toolkit_error}.  

\subsection{Computational cost}
In a second test, we compare the computational costs of the
methods. To do this, we look for the
values of $N=\frac{T}{\Delta t}$ and $m=\frac{\varepsilon_{\max}}{\Delta \varepsilon}$ that enable to reach a
fixed arbitrary error of $Tol=5.10^{-3}$ (recall that in any case the error
cannot exceed 2). For sake of simplicity, we only test powers of $2$.
In this test, we also include the quantified version of the Improved
toolkit II which is described in Remark \ref{quantif}. In our case,
the parameters $\alpha$ and $\beta$ were quantified among 100 values
uniformly distributed in $[0,1]$.

\begin{table}[htbp] 
\begin{center}
\begin{tabular}{|c|c|c|c|}
  \hline
                    & $N=\frac{T}{\Delta t}$&Matrix products  & $m=\frac{\varepsilon_{\max}}{ \Delta \varepsilon}$ \\
  \hline
Strang Op. Splitting&    $16384$            &      32768    &        -                 \\
  \hline
Toolkit             &    $8192$             &      8192     &   $16384$               \\
  \hline
Improved toolkit I  &    $1024$             &      2048     &   $16384$             \\
  \hline
Improved toolkit II &    $4096$             &      12288    &   $16$ \\
  \hline
Quantified Improved toolkit II &    $4096$  &      4096     &   $6400$
\\
  \hline
\end{tabular}
\end{center}
\caption{Values of numerical parameters corresponding to a tolerance
  error of  $Tol=5.10^{-3}$.}\label{tab2} 
\end{table}

These tests show that toolkit methods always give better results as the
second order Strang operator splitting. \\
The two improvements we propose in this paper enable to reduce respectively the global number of
matrix products and the size of the toolkit, which is in agreement
with the analysis we have done. Note that the second improvement reduce
significantly preprocessing step. This fact makes feasible the quantified
version of it, which requires intrinsically a larger toolkit.  

\section*{Acknowledgments}
This work is partially supported by the ANR project C-QUID, INRIA
project ``MicMac'' and by a PICS CNRS-NFS grant.
\bibliographystyle{siam}

\bibliography{references}

\begin{thebibliography}{10}

\bibitem{hbref1}
{\sc A.~Assion, T.~Baumert, M.~Bergt, T.~Brixner, B.~Kiefer, V.~Seyfried,
  M.~Strehle, and G.~Gerber}, {\em Control of chemical reactions by
  feedback-optimized phase-shaped femtosecond laser pulses}, Science, 282
  (1998), pp.~919--922.

\bibitem{kurti-manby-ren-artamonov-ho-rabitz-05}
{\sc G.~Balint-Kurti, F.~Manby, Q.~Ren, M.~Artamonov, T.~Ho, and H.~Rabitz},
  {\em Quantum control of molecular motion including electronic polarization
  effects with a two-stage toolkit}, J. Chem. Phys., 122 (2005).

\bibitem{BST}
{\sc M.~Belhadj, J.~Salomon, and G.~Turinici}, {\em Monotonic time-discretized
  schemes in quantum control{Monotonic time-discretized schemes in quantum
  control}, journal = {J. Phys. A}, year = {2008}, volume = {41}, number = {},
  pages = {362001-362011}, optnote = {DOI: 10.1088/1751-8113/41/36/362001}}.

\bibitem{descombes}
{\sc C.~Besse, B.~Bid{\'e}garay, and S.~Descombes}, {\em Order estimates in
  time of splitting methods for the nonlinear {S}chr\"odinger equation}, SIAM
  J. Numer. Anal., 40 (2002), pp.~26--40 (electronic).

\bibitem{Brezis}
{\sc H.~Brezis}, {\em Analyse fonctionnelle : Thorie et applications}, Dunod,
  1994.

\bibitem{caze}
{\sc T.~Cazenave and A.~Haraux}, {\em An Introduction to Semilinear Evolution
  Equations}, Oxford University, 1998.

\bibitem{hbref90}
{\sc I.~L. Chuang, R.~Laflamme, P.~W. Shor, and W.~H. Zurek}, {\em Quantum
  computers, factoring, and decoherence}, Science, 270 (1995), pp.~1633--1635.

\bibitem{handbook9}
{\sc P.~G. Ciarlet and J.~L. Lions}, eds., {\em Handbook of numerical analysis.
  {V}ol. {IX}}, Handbook of Numerical Analysis, IX, North-Holland, Amsterdam,
  2003.
\newblock Numerical methods for fluids. Part 3.

\bibitem{lubich}
{\sc T.~Jahnke and C.~Lubich}, {\em Error bounds for exponential operator
  splittings}, BIT, 40 (2000), pp.~735--744.

\bibitem{hbref110}
{\sc R.~Judson and H.~Rabitz}, {\em Teaching lasers to control molecules},
  Phys. Rev. Lett, 68 (1992), p.~1500.

\bibitem{hbref3}
{\sc R.~J. Levis, G.~Menkir, and H.~Rabitz}, {\em Selective bond dissociation
  and rearrangement with optimally tailored, strong-field laser pulses},
  Science, 292 (2001), pp.~709--713.

\bibitem{hbref10}
{\sc H.~Rabitz, R.~de~Vivie-Riedle, M.~Motzkus, and K.~Kompa}, {\em Wither the
  future of controlling quantum phenomena?}, Science, 288 (2000), pp.~824--828.

\bibitem{rabitz3}
{\sc H.~Rabitz, G.~Turinici, and E.~Brown}, {\em Control of quantum dynamics:
  Concepts, procedures and future prospects}, in Computational Chemistry,
  Special Volume (C. Le Bris Editor) of Handbook of Numerical Analysis, vol X,
  P.~G. Ciarlet, ed., Elsevier Science B.V., 2003, pp.~833--887.

\bibitem{R-S}
{\sc M.~Reed and B.~Simon}, {\em Methods of Modern Mathematical Physics, II,
  Fourier analysis, self-adjointness}, Academic Press, 1975.

\bibitem{claudedionjcp}
{\sc J.~Salomon, C.~Dion, and G.~Turinici}, {\em Optimal molecular alignment
  and orientation through rotational ladder climbing}, J. Chem. Phys., 123
  (2005), p.~144310.

\bibitem{sanzserna}
{\sc J.~M. Sanz-Serna and J.~G. Verwer}, {\em Conservative and nonconservative
  schemes for the solution of the nonlinear {S}chr\"odinger equation}, IMA J.
  Numer. Anal., 6 (1986), pp.~25--42.

\bibitem{sayood}
{\sc K.~Sayood}, {\em Introduction to data compression, $3^{rd}$ edition},
  Morgan Kaufmann, Elsevier., 2006.

\bibitem{strangsplitting}
{\sc G.~Strang}, {\em On the construction and comparison of difference
  schemes}, SIAM J. Numer. Anal., 5 (1968), pp.~506--517.

\bibitem{hbref9}
{\sc W.~Warren, H.~Rabitz, and M.~Dahleh}, {\em Coherent control of quantum
  dynamics: The dream is alive}, Science, 259 (1993), pp.~1581--1589.

\bibitem{hbref4}
{\sc T.~Weinacht, J.~Ahn, and P.~Bucksbaum}, {\em Controlling the shape of a
  quantum wavefunction}, Nature, 397 (1999), pp.~233--235.

\bibitem{yip-mazziotti-rabitz-03}
{\sc F.~Yip, D.~Mazziotti, and H.~Rabitz}, {\em A local-time algorithm for
  achieving quantum control}, J. Phys. Chem. A., 107 (2003), pp.~7264--7269.

\bibitem{yip-mazziotti-rabitz-03bis}
\leavevmode\vrule height 2pt depth -1.6pt width 23pt, {\em A propagation
  toolkit to design quantum control}, J. Chem. Phys., 118(18) (2003),
  pp.~8168--8172.

\end{thebibliography}

\end{document}